\documentclass[a4paper,12pt]{article}
\usepackage{comment}
\usepackage{cite}
\usepackage{amsmath}
\usepackage{amssymb}
\usepackage{amsfonts}
\usepackage[T1]{fontenc}
\usepackage[utf8]{inputenc}
\usepackage{graphicx}
\usepackage{fancyhdr}
\usepackage{float}
\usepackage{xcolor}
\usepackage{authblk}
\usepackage{mathrsfs}
\usepackage{empheq}
\usepackage[hyphens]{url}
\usepackage{hyperref} 
\usepackage[]{breakurl}

\usepackage{geometry}
\begin{document}

\title{Two integral representations for the logarithm of the Glaisher-Kinkelin constant}

\author[$\dagger$]{Jean-Christophe {\sc Pain}$^{1,2,}$\footnote{jean-christophe.pain@cea.fr}\\
\small
$^1$CEA, DAM, DIF, F-91297 Arpajon, France\\
$^2$Universit\'e Paris-Saclay, CEA, Laboratoire Mati\`ere en Conditions Extr\^emes,\\ 
91680 Bruy\`eres-le-Ch\^atel, France
}

\maketitle

\begin{abstract}
We present two integral representations of the logarithm of the Glaisher-Kinkelin constant. Both are based on a definite integral of $\ln\left[\Gamma(x+1)\right]$, $\Gamma$ being the usual Gamma function. The first one relies on an integral representation of $\ln\left[\Gamma(x+1)\right]$ due to Binet, and the second one results from the so-called Malmst\'en formula. The numerical evaluation is easier with the latter expression than with the former.
\end{abstract}

\section{Introduction}

The Glaisher-Kinkelin constant $A$ is related to the Barnes $G$-function: 
\begin{equation}
    G(n)=\prod_{k=1}^{n-2}k!
\end{equation}
by
\begin{equation}
    A=\lim _{n\rightarrow \infty }{\frac {\left(2\pi \right)^{\frac {n}{2}}n^{{\frac {n^{2}}{2}}-{\frac {1}{12}}}~e^{-{\frac {3n^{2}}{4}}+{\frac {1}{12}}}}{G(n+1)}}.
\end{equation}
The logarithm of the Glaisher-Kinkelin constant $A$ can be expressed as the following integral \cite{Glaisher1878,Almkvist1998}:
\begin{equation}
    \ln A=\frac{1}{12}-2\int_0^{\infty}\frac{x\ln x}{e^{2\pi x}-1}~\mathrm{d}x.
\end{equation}
Another integral is given by \cite{Glaisher1878}:
\begin{equation}\label{gla}
    \int_0^{1/2}\ln\left[\Gamma(x+1)\right]~\mathrm{d}x=-\frac{1}{2}-\frac{7}{24}\ln 2+\frac{1}{4}\ln\pi+\frac{3}{2}\ln A, 	
\end{equation}
where $\Gamma$ represents the usual Gamma function. Only a few integral representations are available in the literature (see for instance the work of Choi and Nash \cite{Choi1997}). Using two different improper integral representations of $\ln\left[\Gamma(x+1)\right]$ (the first one due to Binet (see Sec. \ref{sec1}), and the second one to Malmst\'en  (see Sec. \ref{sec2})), we derive two integral representations of the logarithm of the Glaisher-Kinkelin constant $\ln A$, which, to our knowledge, were not published elsewhere.

\section{Integral representation deduced from the Binet formula}\label{sec1} 

The following important result is due to Binet \cite{Whittaker1990,Sasvari1999}: 
\begin{equation}
    \Gamma(x+1)=\left(\frac{x}{e}\right)^x~\sqrt{2\pi x}~e^{\theta(x)},
\end{equation}
where
\begin{equation}
    \theta(x)=\int_0^{\infty}\left(\frac{1}{e^t-1}-\frac{1}{t}+\frac{1}{2}\right)~\frac{e^{-xt}}{t}~\mathrm{d}t.
\end{equation}
One has
\begin{equation}
    \ln\left[\Gamma(x+1)\right]=x\ln x-x+\frac{1}{2}\ln(2\pi x)+\theta(x),
\end{equation}
and
\begin{equation}
    \int_0^{1/2}\ln\left[\Gamma(x+1)\right]~\mathrm{d}x=\int_0^{1/2}(x\ln x-x)~\mathrm{d}x+\frac{1}{2}\int_0^{1/2}\ln(2\pi x)~\mathrm{d}x+\int_0^{1/2}\theta(x)~\mathrm{d}x.
\end{equation}
The first integral is equal to
\begin{equation}\label{i1}
    \int_0^{1/2}(x\ln x-x)~\mathrm{d}x=-\frac{1}{8}\left(\ln 2+\frac{3}{2}\right),
\end{equation}
the second one is
\begin{equation}\label{i2}
    \int_0^{1/2}\ln(2\pi x)~\mathrm{d}x=\frac{1}{2}\left(\ln \pi-1\right),
\end{equation}
and the third one, applying the Fubini theorem, reads
\begin{equation}\label{i3}
    \int_0^{1/2}\theta(x)~\mathrm{d}x=\int_0^{\infty}\left(\frac{1}{e^t-1}-\frac{1}{t}+\frac{1}{2}\right)~\frac{1}{t}~\left(\int_0^{1/2}\frac{e^{-xt}}{t}~\mathrm{d}x\right)~\mathrm{d}t,
\end{equation}
which is equal to
\begin{equation}\label{i3}
    \int_0^{1/2}\theta(x)~\mathrm{d}x=\int_0^{\infty}\left(\frac{1}{e^t-1}-\frac{1}{t}+\frac{1}{2}\right)~\frac{\left(1-e^{-t/2}\right)}{t^2}~\mathrm{d}t,
\end{equation}
or also
\begin{equation}\label{i3}
    \int_0^{1/2}\theta(x)~\mathrm{d}x=\int_0^{\infty}\frac{\left(1-e^{-t/2}\right)\left[t~\coth\left(t/2\right)-2\right]}{2t^3}~\mathrm{d}t.
\end{equation}
Inserting Eqs. (\ref{i1}), (\ref{i2}) and (\ref{i3}) in Eq. (\ref{gla}) yields
\begin{align}
    -\frac{1}{2}-\frac{7}{24}\ln 2+\frac{1}{4}\ln\pi+\frac{3}{2}\ln A&=-\frac{1}{8}\left(\ln 2+\frac{3}{2}\right)+\frac{1}{4}\left(\ln \pi-1\right)\nonumber\\
    &+\int_0^{\infty}\frac{\left(1-e^{-t/2}\right)\left[t~\coth\left(t/2\right)-2\right]}{2t^3}~\mathrm{d}t.
\end{align}
or equivalently
\begin{align}\label{res1}
    \ln A&=\frac{1}{9}\ln 2+\frac{1}{24}+\frac{2}{3}\int_0^{\infty}\frac{\left(1-e^{-t/2}\right)\left[t~\coth\left(t/2\right)-2\right]}{2t^3}~\mathrm{d}t,
\end{align}
which is the first main result of the present work.

\section{Integral representation resulting from the Malmst\'en formula}\label{sec2}

The Malmst\'en formula reads \cite{Erdelyi1981}:
\begin{equation}
    \ln\Gamma(z+1)=\int_0^{\infty}\left[z-\frac{(1-e^{-zt})}{(1-e^{-t})}\right]\frac{e^{-t}}{t}dt.
\end{equation}
Using the Fubini theorem, one can write
\begin{equation}
    \int_0^{1/2}\ln\left[\Gamma(x+1)\right]~\mathrm{d}x=\int_0^{\infty}\frac{e^{-t}}{t}~\left(\int_0^{1/2}\left[x-\frac{(1-e^{-xt})}{(1-e^{-t})}\right]~\mathrm{d}x\right)~\mathrm{d}t,
\end{equation}
which is equal to
\begin{equation}
    \int_0^{1/2}\ln\left[\Gamma(x+1)\right]~\mathrm{d}x=\int_0^{\infty}\left[\frac{1}{8}-\frac{1}{2\left(1-e^{-t}\right)}+\frac{1}{t\left(1+e^{-t/2}\right)}\right]~\frac{e^{-t}}{t}~\mathrm{d}t,
\end{equation}
or equivalently
\begin{equation}
    \int_0^{1/2}\ln\left[\Gamma(x+1)\right]~\mathrm{d}x=\int_0^{\infty}\frac{e^{-t}\left[(8-3t)~e^t-8~e^{t/2}-t\right]}{8t^2(e^{t}-1)}~\mathrm{d}t.
\end{equation}
Inserting the latter expression in Eq. (\ref{gla}) yields
\begin{equation}
    -\frac{1}{2}-\frac{7}{24}\ln 2+\frac{1}{4}\ln\pi+\frac{3}{2}\ln A=\int_0^{1/2}\ln\left[\Gamma(x+1)\right]~\mathrm{d}x=\int_0^{\infty}\frac{e^{-t}\left[(8-3t)~e^t-8~e^{t/2}-t\right]}{8t^2(e^{t}-1)}~\mathrm{d}t, 	
\end{equation}
leading to
\begin{equation}\label{res2}
    \ln A=\frac{1}{3}+\frac{7}{36}\ln 2-\frac{1}{6}\ln\pi+\frac{2}{3}\int_0^{\infty}\frac{e^{-t}\left[(8-3t)~e^t-8~e^{t/2}-t\right]}{8t^2(e^{t}-1)}~\mathrm{d}t, 	
\end{equation}
which is the second main result of the present work. The convergence of the latter integral (Eq. (\ref{res2})), is much faster than for the representation (\ref{res1}).

\section{Conclusion}

In this letter, we presented two integral representations of the logarithm of the Glaisher-Kinkelin constant. Both are based on a definite integral representation involving $\ln\left[\Gamma(x+1)\right]$, which can itself be expressed by an improper integral, either using a representation proposed a long time ago by Binet, or using the Malmst\'en formula. The two new expressions may help investigating new properties of the Glaisher-Kinkelin constant.


\begin{thebibliography}{99}
    
\bibitem{Glaisher1878} J. W. L. Glaisher, {\it On the product $1^1.2^2.3^3...n^n$}, Messenger Math. {\bf 7}, 43-47 (1878).

\bibitem{Almkvist1998} G. Almkvist, {\it Asymptotic formulas and generalized Dedekind sums}, Experim. Math. {\bf 7}, 343-359 (1998).

\bibitem{Choi1997} J. Choi and C. Nash, {\it Integral representations of the Kinkelin's constant $A$}, Math. Japon. {\bf 45}, 223-230 (1997).

\bibitem{Sasvari1999} Z. Sasvari, {\it An elementary proof of Binet's formula for the Gamma function}, Amer. Math. Mon. {\bf 106}, 156-158 (1999).

\bibitem{Whittaker1990} E. T. Whittaker and G. N. Watson, {\it A course in modern analysis}, 4th ed. Cambridge, England: Cambridge University Press, 1990.

\bibitem{Erdelyi1981} A. Erd\'elyi, W. Magnus, F. Oberhettinger and F. G. Tricomi, {\it Higher transcendental functions}, Vol. 1. New York: Krieger, pp. 20-21, 1981.

\end{thebibliography}
\end{document}